\documentclass{tac}

\usepackage[colorlinks=true]{hyperref}
\hypersetup{allcolors=[rgb]{0.1,0.1,0.4}}

\usepackage{amsmath, mathtools, amssymb}
\usepackage{tikz-cd}

\author{Anthony Bordg}

\thanks{This work was supported by the ERC Advanced Grant ALEXANDRIA (Project GA 742178). I thank Vladimir Voevodsky who commented on an early draft of this document as well as the anonymous referees for their useful comments.}

\address{Department of Computer Science and Technology, University of Cambridge\\
 William Gates Building, JJ Thomson Avenue, Cambridge CB3 0FD, UK
}

\title{The  Interpretation Lifting Theorem for C-Systems}

\copyrightyear{2021}


\eaddress{apdb3@cam.ac.uk}

\newtheorem{theorem}{Theorem}

\newtheoremrm{rem}{Remark}

\newtheoremrm{notation}{Notation}

\newtheoremrm{problem}{Problem}

\newtheoremrm{construction}{Construction}

\newtheoremrm{example}{Example}

\newtheoremrm{conjecture}{Conjecture}

\mathrmdef{Hom}
\mathbfdef{Set}
\mathbfdef{Cat}
\mathbfdef{N}
\mathbfdef{2}
\mathbfdef{CCat}
\mathbfdef{UCat}
\newcommand{\sr}{\rightarrow}

\mathrmdef[cc]{CC}
\mathcaldef[c]{C}
\mathrmdef[preshv]{PreShv}
\mathrmdef[ob]{Ob}
\mathrmdef[mor]{Mor}
\mathrmdef{cod}
\mathrmdef[id]{Id}
\mathrmdef{ft}
\mathrmdef[lan]{Lan}
\mathrmdef{pt}

\mathrmdef{int}
\mathrmdef{op}
\mathrmdef{tt}
\mathrmdef{el}

\begin{document}

\maketitle
\begin{abstract}
In this article we present a solution to a conjecture of Vladimir Voevodsky regarding C-systems. This conjecture provides, under some assumptions, a lift of a functor $M\colon \cc \sr \c$, where $\cc$ is a C-system and $\c$ a category, to a morphism of C-systems $M'\colon \cc \sr \cc(\widehat{\c},p_M)$. We explain the motivation behind this conjecture and introduce the required background material on C-systems. Finally, we give a proof of this conjecture.  
\end{abstract}

\section{Introduction}\label{sec-Introduction}

The late Vladimir Voevodsky devoted the last years of his work to the mathematical theory of type theories. Voevodsky's goal was to give existing type theories a sound mathematical basis that could also apply to future extensions of these type theories. A type theory is a collection of inference rules that can be used as the underlying logic of a proof assistant in order to check mechanically the correctness of mathematical proofs. Occasionally one may want to add a new axiom to this underlying logic. In 2006 Voevodsky proposed to add a new axiom, the \textit{Univalence Axiom}, to the so-called Martin-L\"of type theory~\cite{ml84}. Voevodsky named his new type theory the \textit{Univalent Foundations} (UF) of mathematics~\cite{hottbook}. These new foundations are used for the development of many libraries including \textit{UniMath}~\cite{unimath}, a library of mechanized mathematics in the univalent style using (a version of) the proof assistant Coq based on the Calculus of Inductive Constructions \cite{pm93}, a type theory that is already an extension of Martin-L\"of type theory.  
With the addition of a new axiom such as the Univalence Axiom, one has to prove the soundness of the resulting system. One also wants to give some mathematical interpretations of this system by providing a suitable notion of models for its inference rules. These models are categories equipped with additional operations that correspond to the inference rules of the type theory. Following this approach, one builds a suitable category whose objects are the said models and whose morphisms are functors satisfying some additional properties. Among these models, the model built from the ``raw'' syntax of the type theory is called the \textit{term model}. Central to this approach is the expected result that the term model is an initial object in the category of models. In the case of UF, this expected result is known as the Initiality Conjecture. For a variant of the Calculus of Constructions~\cite{ch86} the corresponding result of initiality was proved by Thomas Streicher in 1988~\cite{streicher91}. Since there are many type theories and a given type theory can be extended to a variety of systems by the addition of new rules, it would be extremely tedious to prove anew the corresponding results of initiality. Instead, Voevodsky wanted to develop a mathematical theory of type theories that would allow to obtain these foundational results ``by specialization of general theorems and constructions for abstract objects the instances of which combine together to produce a given type system''~\cite{vv17}. This program is an instance of building a general theory as a well-motivated problem-solving strategy instead of an ad hoc solution to a given mathematical problem.  

At the heart of Voevodsky's program to achieve this mathematical theory of type theories and prove the Initiality Conjecture lies the notion of a \textit{C-system}, the notion of model developed by Voevodsky. Before being slightly reformulated and developed further by Voevodsky, C-systems were first devised by John Cartmell under the name \textit{contextual categories}~\cite{Cartmell78, Cartmell86} and later studied by Streicher, hence the \textit{C} in \textit{C-system} standing for both Cartmell and \textit{contextual}. The construction of the canonical model of UF in the category of simplicial sets still relies today on the initiality conjecture for contextual categories which remains open~\cite[Conjecture 1.2.9]{kl20}. 
In addition to the Initiality Conjecture for C-systems, Voevodsky formulated another conjecture regarding C-systems in the third article~\cite{vv2015} in his series devoted to this topic~\cite{vv2015identity, vv2016relmonad, vv2016prod, vv2016, vv2017presheaves, vv2017, vf2020}. Unfortunately, this second conjecture was published at the very end of this long article without any explanations~\cite[6.15 Conjecture]{vv2015} and therefore our goal in the first part of this paper consists in giving a more accessible account of this conjecture. For reasons that shall become clear and in order to refer conveniently to this conjecture, we shall name it the \textit{interpretation lifting conjecture}. 

In Section~\ref{sec-conjecture} we shall present the interpretation lifting conjecture and recall the relevant definitions and results in order to put the conjecture in its proper context. This section should make our paper reasonably self-contained. Moreover, no knowledge of the syntax of type theory will be required for the understanding of the conjecture and we will work in set-theoretic foundations as is common in mathematics. Finally, Section~\ref{sec-Solution} will provide a solution to the interpretation lifting conjecture.                          

\section{The Interpretation Lifting Conjecture}\label{sec-conjecture}

In this section we shall present the interpretation lifting conjecture. We start by introducing the relevant background material, proving in the process that one important construction of Voevodsky is actually functorial and that some families of morphisms he introduced are natural transformations. 

\begin{notation}\label{notation-composition}
\begin{enumerate}
\item In order to avoid confusing readers, we will use the standard order for the composition $\circ$ of morphisms, unlike Voevodsky who used the diagrammatic order in his series of papers on C-systems.
\item The category $[{\cal C}^{\op}, \Set]$ of presheaves of sets on $\cal C$ will be denoted $\widehat{\cal C}$.
\end{enumerate} 
\end{notation}

\definition[{C0-system \cite[Definition 2.1]{vv2016}}]\label{def-C0-system}
A C0-system is a category $\cc$ together with the following structure

\begin{enumerate}
\item a function $l \colon \ob(\cc) \sr \mathbb{N}$ named ``length''
\item an object $\pt$ named ``point''
\item a map $\ft \colon \ob(\cc) \sr \ob(\cc)$, with $\ft(X)$ called ``the father of $X$''
\item for each $X \in \ob(\cc)$ a morphism $p_X \colon X \sr \ft(X)$
\item for each $X \in \ob(\cc)$ such that $l(X) > 0$ and each morphism $f \colon Y \sr \ft(X)$ 
an object $f^*X$ and a morphism $q(f,X) \colon f^*X \sr X$
\end{enumerate}
satisfying the following conditions:

\begin{enumerate}
\item $l^{-1}(0)=\{\pt\}$
\item for $X$ such that $l(X) > 0$ one has $l(\ft(X)) = l(X)-1$
\item $\ft(\pt) = \pt$
\item $\pt$ is a final object of the category $\cc$
\item for $X \in \ob(\cc)$ such that $l(X) > 0$ and $f \colon Y \sr \ft(X)$ one has $l(f^*(X)) > 0$, $\ft(f^*X) = Y$ and the distinguished square
\[
\begin{tikzcd}
f^*X \arrow[r, "{q(f,X)}"] \arrow[d, "p_{f^*X}"]
& X \arrow[d, "p_X"] \\
Y \arrow[r, "f"] 
& \ft(X)
\end{tikzcd}
\]
commutes
\item for $X \in \ob(\cc)$ such that $l(X) > 0$ one has $(\id_{\ft(X)})^*X = X$ and $q(\id_{\ft(X)},X) = \id_X$
\item for $X \in \ob(\cc)$ such that $l(X) > 0$, $g \colon Z \sr Y$ and $f \colon Y \sr \ft(X)$ one has $(f \circ g)^*X = g^*(f^*X)$ and $q(f \circ g,X) = q(f,X) \circ q(g,f^*X)$
\end{enumerate}
\enddefinition

For every morphism $f \colon Y \sr X$, the morphism $\ft(f) \colon Y \sr \ft(X)$ will denote the post-composition of $f$ with $p_X$. 

\definition[C-system {\cite[Definition 2.3]{vv2016}}]\label{def-C-system}
A C-system is a C0-system equipped with an operation $f \mapsto s_f$ defined for all $f \colon Y \sr X$ such that $l(X) > 0$ and satisfying the following properties.

\begin{enumerate}
\item $s_f \colon Y \sr (\ft(f))^*X$
\item $p_{(\ft(f))^*X} \circ s_f = \id_Y$
\item $q(\ft(f),X) \circ s_f = f$
\item if $X = g^*U$, where $g \colon \ft(X) \sr \ft(U)$, then $s_{q(g,U) \circ f} = s_f$
\end{enumerate}
The map $s_f$ will be called the section of $f$.
\enddefinition

The reader can check that every distinguished square in a C-system is a pullback square. It is actually equivalent for a C0-system $\cc$ to be a C-system and for its distinguished squares to be pullback squares \cite[Proposition 2.4]{vv2016}. 

\example[{\cite{vv1bonn2016}}]\label{example-C-systems}
Let $\N_{triv}$ be the category with set of objects the set $\mathbb{N}$ of natural numbers and with exactly one morphism between any two objects. There exists a C-system structure on $\N_{triv}$ given by the identity map as the length function. The other operations are then completely determined.  
\endexample

\begin{remark}\label{remark-equivalence-of-categories}
The reader should note that C-system structures cannot be transported along equivalences of categories. Indeed, consider the category $\2_{triv}$ with two objects and one isomorphism between them. This category is equivalent to the category $\N_{triv}$, but the reader can check that there does not exist a C-system structure on $\2_{triv}$ \cite{vv1bonn2016}. C-system structures being algebraic structures, the right notion of sameness for C-systems is the notion of an isomorphism. 
\end{remark}

\definition[{\cite[Remark 2.8]{vv2016}}]\label{morphism-of-C-systems}
Let $(\cc,l,\pt,\ft,p,q,s)$ and $(\cc',l',\pt',\ft',p',q',s')$ be two C-systems, a morphism of C-systems is a functor $F: \cc \sr \cc'$ that respects the length functions, the final objects, the $p$-operations, the $s$-operations and commutes with the father functions and the $q$-operations whenever these functions and operations are defined. In other words the following equalities

\begin{align*}
l'(F(X)) &= l(X) \\
F(p_X) &= p_{F(X)} \\
F(s_f) &= s_{F(f)} \\
\ft'(F(X)) &= F(\ft(X)) \\
q'(F(f),F(X)) &= F(q(f,X))
\end{align*} 
are satisfied.
\enddefinition

\begin{remark}
Such a functor $F$ automatically satisfies $F(\pt) = \pt'$. 
\end{remark}

The category of C-systems and their morphisms will be denoted $\CCat$.
We shall now introduce the notion of a \textit{universe category} that will play an important role in the next section.

\definition[universe category {\cite[2.6 Definition]{vv2015}}]
A universe category is a triple $(\mathcal{C}, p, \pt)$, often denoted simply by $(\mathcal C, p)$, where $\cal C$ is a category, $\pt$ is a final object in $\cal C$ and $p \colon \widetilde{U} \sr U$ is a morphism in $\cal C$ together with, for every morphism $f \colon X \sr U$, a chosen pullback square as follows.

\[
\begin{tikzcd}
(X;f) \arrow[r, "Q(f)"] \arrow[d, "p_{X,f}"] 
& \widetilde{U} \arrow[d, "p"] \\
X \arrow[r, "f"] 
& U
\end{tikzcd}
\] 
\enddefinition 

\begin{example}
Consider $(\N_{triv})^{\op}$ the opposite category of $\N_{triv}$ (\textit{cf.} \ref{example-C-systems}), $p\colon 1 \sr 0$ the unique morphism from $1$ to $0$ and for every morphism $f\colon n \sr 0$ take $(n;f)\coloneqq n+1$, then $Q(f)$ (\textit{resp.} $p_{n,f}$) is the unique morphism from $n+1$ to $1$ (\textit{resp.} from $n+1$ to $n$). Note that every morphism in $(\N_{triv})^{\op}$ is an isomorphism. The triple $((\N_{triv})^{\op}, p, 0)$ is a universe category, since a commutative square where all four arrows are isomorphisms is a pullback square. 
\end{example}  

Voevodsky proved that one can define a C-system from a universe category $(\mathcal{C}, p, \pt)$~\cite[2.12 Construction]{vv2015}. This C-system will be denoted $\cc(\mathcal{C}, p, \pt)$ and following Voevodsky it will often be abbreviated to $\cc(\mathcal{C}, p)$. 
Moreover, the C-system $\cc(\mathcal{C}, p)$ comes equipped with a fully faithful functor $\int\colon \cc(\mathcal{C}, p) \sr \mathcal{C}$ from the underlying category of $\cc(\mathcal{C}, p)$ to $\mathcal{C}$~\cite[2.9 Lemma]{vv2015}. For the convenience of the reader, we shall briefly recapitulate these constructions. 

\begin{construction}\label{construction-of-generated-C-system}
We first define sets $\ob_n(\mathcal{C},p)$, shortened $\ob_n$, and maps 
\[
\int_n\colon \ob_n \sr \ob(\mathcal{C})\] 
by a mutual recursion. The set $\ob_0\coloneqq \{\tt\}$ is a distinguished singleton with $\tt$ its unique element, $\int_0$ maps the unique element of $\ob_0$ to $\pt$ and the recursive cases are given as follows
\[
\ob_{n+1}\coloneqq \coprod_{A \in \ob_n} \Hom_{\mathcal{C}}(\int_n(A), U) 
\]
and
\[
\int_{n+1}(A, f) \coloneqq (\int_n(A); f).
\]
The set $\ob(\cc(\mathcal{C},p))$ of objects of $\cc(\mathcal{C},p)$ is then simply 
\[
\displaystyle \coprod_n \ob_n
\]
with the length function being the obvious projection and $\int$ on objects being the sum of the maps $\int_n$, while the set of morphisms $\mor(\cc(\mathcal{C},p))$ of $\cc(\mathcal{C},p)$ is 
\[
\coprod_{\Gamma,\Gamma' \in \ob(\cc(\mathcal{C},p))} \Hom_{\mathcal{C}}(\int(\Gamma), \int(\Gamma')),
\]
the functor $\int$ mapping a morphism $(\Gamma, \Gamma', a)$ to $a$. The point of $\cc(\mathcal{C},p)$ is $(0,\tt)$. The father function $\ft$ is the sum of the maps $\ft_n$, where $\ft_0 \coloneqq \id_{\ob_0}$ and $\ft_{n+1}$ maps an object $(A,f)$ of $\ob_{n+1}$ to $A$ in $\ob_n$. Last, we have to define the distinguished pullback squares of the C-system $\cc(\mathcal{C},p)$. First, we need a morphism $p_{\Gamma}\colon \Gamma \sr \ft(\Gamma)$ for every $\Gamma \in \ob(\cc(\mathcal{C},p))$. For $\Gamma \in \ob_0$, take $p_{\Gamma} \coloneqq \id_{\Gamma}$ and for $\Gamma \coloneqq (n+1, A)$ in $\ob_{n+1}$ with $A \coloneqq (B, f)$, take $p_{\Gamma}\coloneqq (\Gamma, \ft(\Gamma), p_{\int_n(B),f})$. Second, for each $\Gamma \in \ob(\cc(\mathcal{C},p))$ such that $l(\Gamma) > 0$ and each morphism $(\Gamma', \ft(\Gamma), f)$, we need an object $f^*\Gamma$ and a morphism $q(f,\Gamma)\colon f^*\Gamma \sr \Gamma$. Assume $\Gamma$ is $(n+1, A)$, with $A \coloneqq (B, g)$ in $\ob_{n+1}$, and assume $\Gamma'$ is $(m, C)$. In this case $g$ is a morphism in $\mathcal{C}$ from $\int_n(B)$ to $U$, while $f$ is a morphism from $\int_m(C)$ to $\int_n(B)$. Take $f^*\Gamma \coloneqq (m+1, (C, g\circ f))$ and $q(f, \Gamma)$, seen as an arrow in $\mathcal{C}$, is the dashed arrow obtained from the universal property of the pullback square in the following diagram.
\[
\begin{tikzcd}[column sep = large, row sep = large]
(\int_m(C); g\circ f) \arrow[rr, "Q(g\circ f)"] \arrow[rd, dashrightarrow] \arrow[dd, "p_{\int(C), g\circ f}"'] & & \widetilde{U} \arrow[dd, "p"]\\
	& (\int_n(B); g) \arrow[ru, "Q(g)"] \arrow[d, "p_{\int(B),g}"]&	\\
\int_m(C) \arrow[r, "f"]& \int_n(B) \arrow[r, "g"] & U
\end{tikzcd}         
\]
\end{construction}

We are now in a position to state the Interpretation Lifting Conjecture.
 
\conjecture[{\cite[6.15 Conjecture]{vv2015}}]\label{interpretation-lifting-conjecture}
Let $\cal C$ be a category, $\cc$ be a C-system and $M\colon \cc \sr {\cal C}$ a functor such that $M(\pt)$ is a final object of $\cal C$ and $M$ maps the distinguished pullback squares of $\cc$ to pullback squares of $\cal C$. Then there exists a universe category $(\widehat{\cal C}, p_M)$ and a C-system morphism $M'\colon \cc \sr \cc(\widehat{\cal C}, p_M)$ such that the square
\[
\begin{tikzcd}
\cc \arrow[r, "M"] \arrow[d, "M'"] 
& {\cal C} \arrow[d, "Y_{\mathcal{C}}"] \\
\cc(\widehat{\cal C}, p_M) \arrow[r, "\int"]
& \widehat{\mathcal{C}}
\end{tikzcd}
\]
where $Y_{\mathcal{C}}$ is the Yoneda embedding, commutes up to a functor isomorphism. 
\endconjecture
At this point we shall offer a few words of motivation from Voevodsky:
\begin{quote}
Suppose $\cc$ is the syntactic C-system of a type theory. Then a functor such as $M$ is a ``weak interpretation'' of the type theory, because by passing from a C-system that is a rigid algebraic structure defined up to an isomorphism, to a category $\mathcal{C}$ that is a much less rigid structure defined up to an equivalence, we can ``erase'' a lot of structure that exists in $\cc$.
By constructing $M'$ one lifts a ``weak" interpretation to a ``strong" one, with values in a C-system [of the form] $\cc(\mathcal{C}, p)$. Such an interpretation is ``strong" because it respects all the structures of the C-system $\cc$ that are erased by the original functor $M$.\footnote{private communication}
\end{quote}
The reader should note that the ``syntactic C-system of a type theory'' is just another way to refer to what we called in the introduction the \textit{term model} of a type theory which is expected to be an initial object in $\CCat$ (\textit{cf.} Section~\ref{sec-Introduction}). Voevodsky's comment echoes the Remark~\ref{remark-equivalence-of-categories} emphasizing that C-system structures cannot be transported along equivalences of categories.   

A couple of propositions are in order as well as a couple of lemmas that will be useful later in Section \ref{sec-Solution}.
First, note that every C-system $\cc$ is actually the C-system defined from some universe category.

\proposition[{\cite[5.2 Construction]{vv2015}}]\label{prop-iso-of-C-systems}
For every C-system $\cc$, there exists a universe category $(\widehat{\cc}, \partial)$ such that $\cc$ and $\cc(\widehat{\cc}, \partial)$ are isomorphic as C-systems.
\endproposition

We should recall here some details about the universe category $(\widehat{\cc}, \partial)$ for which there exists an isomorphism $I_{\cc}\colon \cc \sr \cc(\widehat{\cc}, \partial)$. 

\begin{construction}\label{construction-of-iso}
Let $U$ be the presheaf that maps an object $\Gamma$ of $\cc$ to the set
\[
\left\lbrace \Delta \mid l(\Delta) > 0 \text{ and } \ft(\Delta) = \Gamma \right\rbrace 
\]
and maps a morphism $f$ to the function $U(f)$ defined by $U(f)(\Delta) \coloneqq f^*\Delta$. Let $\widetilde{U}$ be the presheaf that maps an object $\Gamma$ of $\cc$ to the set
\[
\left\lbrace s \in \mor(\cc) \mid s\colon \ft(\Delta) \sr \Delta,\, l(\Delta) > 0,\, \ft(\Delta) = \Gamma \text{ and } p_{\Delta} \circ s = \id_\Gamma \right\rbrace
\]
of sections of the canonical projections $p_{\Delta}$ for $\Delta$ such that $l(\Delta) > 0$ and $\ft(\Delta) = \Gamma$
and such that $\widetilde{U}$ maps a morphism $f$ to the function $\widetilde{U}(f)$ defined by $\widetilde{U}(f)(s) \coloneqq q(f,\Delta)^*s$. The natural transformation $\partial$ simply maps a section to its codomain. Let $\pt$ be the constant presheaf given by a distinguished singleton $\lbrace \star \rbrace$ in $\Set$. Then $(\widehat{\cc}, \partial, \pt)$ together with the canonical pullback squares in the presheaf category $\widehat{\cc}$ is a universe category. We will construct the isomorphism $I_{\cc}\colon \cc \sr \cc(\widehat{\cc}, \partial)$ as follows. For every $\Gamma$ in $\cc$, the canonical bijection
\[
U(\Gamma) \cong \Hom_{\widehat{\cc}}(Y_{\Gamma}, U)
\] 
given by the Yoneda lemma will be denoted $u_{\Gamma}$. 
Let us denote $\delta(\Delta)$ the section of $p_{\Delta, p_{\Delta}}$ given by the diagonal, its image under the canonical bijection $\widetilde{U}(\Delta) \cong \Hom(Y_{\Delta}, \widetilde{U})$ will be denoted $\widetilde{u}_{\Delta}(\delta(\Delta))$. For every $\Gamma$ in $\cc$ and every $\Delta$ in $U(\Gamma)$,
\[
\gamma_{\Delta}\colon (Y_{\Gamma}; u_{\Gamma}(\Delta)) \sr Y_{\Delta}
\]
will denote the isomorphism given by the universal property of the following pullback square.
\[
\begin{tikzcd}
Y_{\Delta} \arrow[rrd, bend left, "\widetilde{u}_{\Delta}(\delta(\Delta))"] \arrow[rdd, bend right, "p_{\Delta} \circ -"'] \arrow[rd, dashrightarrow] & & \\
	& (Y_{\Gamma}; u_{\Gamma}(\Delta)) \arrow[d] \arrow[r]
	& \widetilde{U} \arrow[d, "\partial"] \\
	& Y_{\Gamma} \arrow[r, "u_{\Gamma}(\Delta)"] & U
\end{tikzcd}
\]
Finally, let us denote $\ob_n(\cc)$ the set of objects in $\cc$ of length $n$. We define pairs $(I_n, \psi_n)$ by a mutual recursion, where $I_n\colon \ob_n(\cc) \sr \ob_n(\widehat{\cc},\partial)$ is a function and $\psi_n(\Gamma)\colon \int(I_n(\Gamma)) \sr Y_{\Gamma}$ is an isomorphism for every $\Gamma$ in $\ob_n(\cc)$. We take $I_0(\pt) = \pt$ and $\psi_0(\pt)$ is the unique isomorphism from our choice of final object $\pt$ in $\widehat{\cc}$ to $Y_{\pt}$. The recursion step is then given for every $\Delta \in U(\Gamma)$ with $I_n(\Gamma) = B$ by the equalities
\begin{align*}
I_{n+1}(\Delta) &= (B, u_{\Gamma}(\Delta) \circ \psi_n(\Gamma)) \\
\psi_{n+1}(\Delta) &= \gamma_{\Delta} \circ Q(\psi_n(\Gamma), u_{\Gamma}(\Delta)),
\end{align*} 
where $Q(\psi_n(\Gamma), u_{\Gamma}(\Delta))$ denotes the dashed arrow obtained from the universal property of the pullback square in the following diagram.
\[
\begin{tikzcd}[column sep = huge, row sep = large]
(\int(I(\Gamma)); u_{\Gamma}(\Delta) \circ \psi_n(\Gamma)) \arrow[r, dashrightarrow]  \arrow[rr, bend left, "Q(u_{\Gamma}(\Delta) \circ \psi_n(\Gamma))"] \arrow[d]
& (Y_{\Gamma}; u_{\Gamma}(\Delta)) \arrow[r, "Q(u_{\Gamma}(\Delta))"] \arrow[d]
& \widetilde{U} \arrow[d, "\partial"] \\
\int(I(\Gamma)) \arrow[r, "\psi_n(\Gamma)"] 
& Y_{\Gamma} \arrow[r, "u_{\Gamma}(\Delta)"]& U
\end{tikzcd}
\]
The isomorphism $I_{\cc}$ maps an object $\Gamma$ to $(l(\Gamma), I_{l(\Gamma)}(\Gamma))$ and a morphism $f\colon \Gamma' \sr \Gamma$ to $(I_{\cc}(\Gamma'), I_{\cc}(\Gamma), \psi(\Gamma)^{-1} \circ Y_{\cc}(f) \circ \psi(\Gamma'))$, where $Y_{\cc}\colon \cc \sr \widehat{\cc}$ denotes the Yoneda embedding. 
\end{construction}

\begin{lemma}\label{lemma-natural-iso-with-yoneda}
There exists a natural isomorphism $\psi$ from $\int \circ I_{\cc}$ to $Y_{\cc}$.
\end{lemma}
\proof
For each object $\Gamma$ of $\cc$, we define a morphism $\psi_{\Gamma}\colon \int(I_{\cc}(\Gamma)) \sr Y_{\Gamma}$ as $\psi_{\Gamma}\coloneqq \psi_{l(\Gamma)}(\Gamma)$ (see Construction \ref{construction-of-iso}). For every morphism $f\colon \Gamma' \sr \Gamma$, we need to prove that the following diagram commutes.
\[
\begin{tikzcd}
\int(I_{\cc}(\Gamma')) \arrow[r, "\psi_{\Gamma'}"] \arrow[d, "\int(I_{\cc}(f))"']
& Y_{\Gamma'} \arrow[d, "Y_{\cc}(f)"] \\
\int(I_{\cc}(\Gamma)) \arrow[r, "\psi_{\Gamma}"] & Y_{\Gamma}
\end{tikzcd}
\]
It is easily checked as follows.
$$\mld \psi_{\Gamma} \circ \int(I_{\cc}(f)) &= \psi_{l(\Gamma)}(\Gamma) \circ \int(I_{\cc}(f)) \\
= \psi_{l(\Gamma)}(\Gamma) \circ \psi_{l(\Gamma)}(\Gamma)^{-1} \circ Y_{\cc}(f) \circ \psi_{l(\Gamma')}(\Gamma') \\
= Y_{\cc}(f) \circ \psi_{l(\Gamma')}(\Gamma') \\
= Y_{\cc}(f) \circ \psi_{\Gamma'}$$
\endproof

We shall define the notion of a morphism of universe categories, which Voevodsky called a \textit{functor of universe categories}~\cite[4.1 Definition]{vv2015}.

\definition[{\cite[4.1 Definition]{vv2015}}]\label{def-morphism-universe-categories}
A morphism between universe categories $(\mathcal{C}, p, \pt)$ and $(\mathcal{C'}, p', \pt')$ is a triple $(F, \phi, \widetilde{\phi})$, where $F\colon \mathcal{C} \sr \mathcal{C'}$ is a functor, $\phi\colon F(U) \sr U'$ and $\widetilde{\phi}\colon F(\widetilde{U}) \sr \widetilde{U'}$ are morphisms in $\mathcal{C'}$, such that $F$ maps the chosen pullback squares based on $p$ to pullback squares, $F(\pt)$ is a final object of $\mathcal{C'}$ and the following square

\[
\begin{tikzcd}
F(\widetilde{U}) \arrow[r, "\widetilde{\phi} "] \arrow[d, "F(p)"]
& \widetilde{U'} \arrow[d, "p'"] \\
F(U) \arrow[r, "\phi"] 
& U'
\end{tikzcd}
\]
is a pullback square. 
\enddefinition

Given two morphisms of universe categories
\[
(F, \phi, \widetilde{\phi})\colon (\mathcal{C}, p, \pt) \sr (\mathcal{C'},p',\pt')
\]
 and 
 \[
 (G, \psi, \widetilde{\psi})\colon (\mathcal{C}',p',\pt') \sr (\mathcal{C''}, p'', \pt''),
 \]
we define their composition as $(G\circ F, \psi \circ G(\phi), \widetilde{\psi} \circ G(\widetilde{\phi}))$. Since two pullback squares based on the same diagram are connected by an isomorphism and given that a functor maps an isomorphism to an isomorphism, one readily checks that the triple 
\[
(G\circ F, \psi \circ G(\phi), \widetilde{\psi} \circ G(\widetilde{\phi}))
\]
is a morphism of universe categories from $(\mathcal{C}, p, \pt)$ to $(\mathcal{C''}, p'', \pt'')$. We define the identity morphism of $(\mathcal{C}, p, \pt)$ as $(\id_{\mathcal{C}}, \id_U, \id_{\widetilde{U}})$. The associativity and unitality of this composition are straightforward. 
The category of universe categories will be denoted $\UCat$.  
Also, from a morphism of universe categories $(F, \phi, \widetilde{\phi})\colon (\mathcal{C}, p, \pt) \sr (\mathcal{C'}, p', \pt')$, it is possible to define a C-system morphism $\cc(F, \phi, \widetilde{\phi})\colon \cc(\mathcal{C},p) \sr \cc(\mathcal{C'},p')$ between the corresponding C-systems (see \cite[4.7 Construction]{vv2015}, where this last morphism is denoted $H$). We shall also recapitulate briefly this construction for the convenience of the reader. 

\begin{construction}\label{construction-lift-of-functor-of-univ-cat}
Let us denote $\psi$ the isomorphism from $\pt'$ to $F(\pt)$. We first define by a mutual recursion maps $H_n\colon \ob_n \sr \ob'_n$ and isomorphisms $\psi_n(A)\colon \int'(H_n(A)) \sr F(\int(A))$ for every $A \in \ob_n$. Take $H_0$ to be the unique map from $\ob_0$ to $\ob'_0$ and $\psi_0(A) \coloneqq \psi$. The recursive cases are given as follows: 
\[
H_{n+1} \coloneqq (H_n(A), \phi \circ F(f) \circ \psi_n(A))
\]
and
\[
\psi_{n+1}(A,f)\colon (\int(H_n(A)); \phi \circ F(f) \circ \psi_n(A)) \sr F(\int(A,f)) 
\]
is the unique morphism in the following diagram
\[
\begin{tikzcd}[column sep = huge, row sep = huge]
\int'(H_{n+1}(A,f)) \arrow[r, "{\psi_{n+1}(A,f)}", dashrightarrow] \arrow[d, "{p_{\int'(H_n(A)), \phi \circ F(f) \circ \psi_n(A)}}"] \arrow[rrr, bend left = 22, "Q(\phi \circ F(f) \circ \psi_n(A))"]
& F(\int(A,f))  \arrow[r, "F(Q(f))"] \arrow[d, "F(p_{\int(A),f})"]
& F(\widetilde{U}) \arrow[r, "\widetilde{\phi}"] \arrow[d, "F(p)"]
& \widetilde{U'} \arrow[d, "p"] \\
\int'(H_n(A)) \arrow[r, "\psi_n(A)"]
& F(\int(A)) \arrow[r, "F(f)"]
& F(U) \arrow[r, "\phi"]
& U'
\end{tikzcd}
\]      
such that the equalities
\begin{align*}
F(p_{\int(A),f}) \circ \psi_{n+1}(A,f) &= \psi_n(A) \circ p_{\int'(H_n(A)), \phi \circ F(f) \circ \psi_n(A)} \\
\widetilde{\phi} \circ F(Q(f)) \circ \psi_{n+1}(A, f) &= Q(\phi \circ F(f) \circ \psi_n(A)) 
\end{align*}
hold. The functor $\cc(F, \phi, \widetilde{\phi}) \coloneqq H$ is then given on objects by the sum of the functions $H_n$, while  on morphism $H$ maps $(\Gamma, \Gamma', f)$ to $(H(\Gamma), H(\Gamma'), \psi(\Gamma')^{-1} \circ F(f) \circ \psi(\Gamma))$.
\end{construction}

\begin{lemma}\label{lemma-natural-iso-with-functor-of-univ-cat}
There exists a natural isomorphism $\psi$ from $\int' \circ H$ to $F \circ \int$.
\end{lemma}
\proof
For each element $A$ of $\ob_n(\mathcal{C},p)$, define $\psi_A$ the component of $\psi$ at $A$ as $\psi_{l(A)}(A)$ (see Construction \ref{construction-lift-of-functor-of-univ-cat}). The argument to show that $\psi$ is a natural transformation is similar to the one in Lemma \ref{lemma-natural-iso-with-yoneda}.  
\endproof
 
\begin{proposition}
The maps $(\mathcal{C}, p, \pt) \mapsto \cc(\mathcal{C}, p, \pt)$ and $(F, \phi, \widetilde{\phi}) \mapsto \cc(F, \phi, \widetilde{\phi})$ define a functor $\cc(-,-,-)$ from $\UCat$ to $\CCat$. 
\end{proposition} 
\proof
We have to prove the equality 
\[
\cc(\id_{\mathcal{C}}, \id_U, \id_{\widetilde{U}}) = \id_{\cc(\mathcal{C},p)},
\]
for every universe category $(\mathcal{C}, p\colon \widetilde{U} \sr U, \pt)$. We have also to prove the equality 
\[
\cc(G\circ F, \psi \circ G(\phi), \widetilde{\psi} \circ G(\widetilde{\phi})) = \cc(G, \psi, \widetilde{\psi}) \circ \cc(F, \phi, \widetilde{\phi}),
\]
namely that one obtains the same morphism of C-systems if one starts by lifting the two morphisms of universe categories and then composes the resulting morphisms of C-systems or if one starts by composing the two morphisms of universe categories and then lifts the resulting morphism of universe categories. Both equalities follow from a proof by induction on $n$ in the formulas defining $H_n$ and $\psi_n$ above.
\endproof

\section{Solution}\label{sec-Solution}

\subsection{Universe categories and left Kan extensions}
\label{subsec-universe-cat}

Let $\cc$ be a C-system, $\mathcal{C}$ a category and $M\colon \cc \sr \mathcal{C}$ a functor from the underlying category of $\cc$ to $\mathcal{C}$ such that $M(\pt)$ is a final object of $\cal C$ and $M$ maps the distinguished pullback squares of $\cc$ to pullback squares of $\cal C$.
Let $(\widehat{\cc}, \partial, \pt)$ be the universe category of Construction \ref{construction-of-iso} together with its isomorphism $I_{\cc}\colon \cc \sr \cc(\widehat{\cc}, \partial)$.

\begin{problem}
To construct a universe category $(\widehat{\cal C},\partial',\pt')$ and a functor of universe categories from $(\widehat{\cc}, \partial, \pt)$ to $(\widehat{\cal C}, \partial', \pt')$.
\end{problem}

\begin{construction}\label{construction-univ-cat}
Consider the functor $M_! \coloneqq \lan_{Y_{\cc}} (Y_{\cal C} \circ M)$ from $\widehat{\cc}$ to $\widehat{\cal C}$, where $\lan_{Y_{\cc}} (Y_{\cal C} \circ M)$ denotes the left Kan extension of $Y_{\cal C} \circ M$ along the (covariant) Yoneda embedding $Y_{\cc}\colon \cc \sr \widehat{\cc}$. We define $\partial'$ as $M_!(\partial)$. Let $(\widehat{\cal C}, \partial')$ be the universe category where the pullback squares based on $\partial'$ are the canonical pullback squares in the presheaf category $\widehat{\cal C}$.  
\end{construction}

\begin{lemma}\label{lemma-preservation-of-final-object}
The object $M_!(\pt)$ is final in $\widehat{\cal C}$. 
\end{lemma}
\proof
Since $M(\pt)$ is a final object by assumption, then $Y_{\cal C}(M(\pt))$ is a final object and the slice category $Y_{\cc}/\pt$ is isomorphic to $\cc$, hence the left Kan extension $M_!$ at $\pt$ is given by the following colimit.
\[
M_!(\pt) = \varinjlim_{x\in \cc} Y_{\cal C}(M(x))
\]
The object $Y_{\cal C}(M(\pt))$ being final, we have an isomorphism
\[
\varinjlim_{x\in \cc} Y_{\cal C}(M(x)) \cong Y_{\cal C}(M(\pt))
\]
so we conclude.
\endproof

Given $y\in \cc$ and $u \in U(y)$, let $\delta(u)$ denote the section obtained from the universal property of the following distinguished pullback square in $\cc$.
\[
\begin{tikzcd}[column sep=large, row sep=large]
u \arrow[rd, dashrightarrow, "\delta(u)"] \arrow[rdd, bend right, "\id"] \arrow[rrd, bend left, "\id"] & & \\
	& {p_{u}}^* u \arrow[r, "{q(p_u,u)}"] \arrow[d, "p_{{p_u}^*u}"] &
	u \arrow[d, "p_u"] \\
	& u \arrow[r, "p_u"] &
y
\end{tikzcd}
\]

\begin{lemma}\label{lemma-pullback-of-diagonal-section}
We have the equality $\widetilde{U}(q(f,u) \circ s) \delta(u) = s$ for every object $x$ of $\cc$, $f\colon x \sr y$ and every $s \in \widetilde{U}(x)$ such that $\partial_x(s) = f^*u$.
\end{lemma}
\proof
By definition of $\widetilde{U}$, the morphism $\widetilde{U}(q(f,u) \circ s) \delta(u)$ is $q(q(f,u) \circ s, {p_u}^*u)^* \delta(u)$, namely the pullback of $\delta(u)$ along the morphism $q(q(f,u) \circ s, {p_u}^*u)$. Since $\partial_x(s) = f^*u$, $s$ is a section of $p_{f^*u}$ and we have the equalities (\textit{cf.} point 7 of Definition \ref{def-C0-system})
$$\mld (q(f,u) \circ s)^* {p_u}^*u &= (p_u \circ q(f,u) \circ s)^* u \\ = f^* u.$$
It means that $\widetilde{U}(q(f,u) \circ s) \delta(u)$ is the unique section $\alpha$ of $p_{f^*u}$ satisfying 
\[
q(q(f,u) \circ s, {p_u}^*u) \circ \alpha = \delta(u) \circ q(f,u) \circ s,\]
hence by unicity it suffices to prove that the equality
\[
q(q(f,u) \circ s, {p_u}^*u) \circ s = \delta(u) \circ q(f,u) \circ s
\]
holds. Consider the following universal problem
\[
\begin{tikzcd}[column sep= huge, row sep=huge]
f^* u \arrow[d, "p_{f^*u}"] \arrow[rr, "{q(q(f,u)\circ s, {p_u}^*u)}"]&
 & {p_u}^* u \arrow[r, "{q(p_u, u)}"] \arrow[d, "p_{{p_u}^*u}"]
 & u \arrow[d, "p_u"]  \\
x \arrow[u, bend left, "s"] \arrow[r, "s"] \arrow[rru, dashrightarrow, "\beta"]& 
f^*u \arrow[r, "{q(f,u)}"]& 
u \arrow[r, "p_u"] \arrow[u, bend left, "\delta(u)"]& 
y,
\end{tikzcd}
\]
where $\beta$ is the unique morphism satisfying the equations
\begin{align*}
p_{{p_u}^* u} \circ \beta &= q(f,u) \circ s \\
q(p_u, u) \circ \beta &= q(p_u,u) \circ q(q(f,u) \circ s, {p_u}^*u) \circ s.
\end{align*}
Since we have the equalities
$$\mld q(p_u,u) \circ q(q(f,u) \circ s, {p_u}^*u)  &= q(p_u \circ q(f,u) \circ s, u) \\ = q(f, u),$$
it is easy to check that both $q(q(f,u) \circ s, {p_u}^*u) \circ s$ and $\delta(u) \circ q(f,u) \circ s$ are solutions of this universal problem, hence they are equal.   
\endproof

Write $U$ as a colimit of representables 
\[
\varinjlim_{(y,u)\in \el(U)^{\op}} y,
\]
where $y$ stands for the representable $Y_{\cc}(y)$ and let $c_{(y,u)}$ denote the edge from the copy of $y$ indexed by $(y,u)$ to $U$ given by the cocone of the latter.
 
\begin{lemma}\label{lemma-pullback-of-univ-along-cocone}
The square
\[
\begin{tikzcd}[column sep=large, row sep=large]
u \arrow[r, "\delta(u)"] \arrow[d, "Y(p_u)"] & 
\widetilde{U} \arrow[d, "\partial"] \\
y \arrow[r, "c_{(y,u)}"] & 
U,
\end{tikzcd}
\] 
where $\delta(u)$ denotes the natural transformation that corresponds to $\delta(u)$ in $\widetilde{U}(u)$, is a pullback square. 
\end{lemma}
\proof
Since limits are pointwise, it suffices to prove that the square 
\[
\begin{tikzcd}[column sep= large]
\Hom(x,u) \arrow[d, "p_u \circ -"] \arrow[r, "{\delta(u)}_x"]& 
\widetilde{U}(x) \arrow[d, "\partial_x"] \\
\Hom(x,y) \arrow[r, "{c_{(y,u)}}_x"] & 
U(x)
\end{tikzcd}
\]
is a pullback square in $\Set$ for every object $x$ of $\cc$. Let 
\[
\phi\colon \Hom(x,y) \times_{U(x)} \widetilde{U}(x) \rightarrow \Hom(x,u)
\] 
be the map sending $(f,s)$, such that $\partial_x(s) = f^*u$, to $q(f,u) \circ s$. Let $\psi$ be the map that sends $g$ to $(p_u \circ g, \widetilde{U}(g)\delta(u))$, where, by definition of $\widetilde{U}$, $\widetilde{U}(g) \delta(u)$ is the pullback of $\delta(u)$ along the morphism $q(g, {p_u}^*u)$. Since $g^*({p_u}^*u)$ is equal to $(p_u \circ g)^*u$ for every $g$ in $\Hom(x,u)$ (by point 7 in Definition \ref{def-C0-system}), the map $\psi$ has values in $\Hom(x,y)\times_{U(x)} \widetilde{U}(x)$. Using $q(p_u \circ g,u) = q(p_u, u) \circ q(g, {p_u}^*u)$ (\textit{cf. ibid}), we conclude $\phi \circ \psi = \id$. Since for every object $(f,s)$ of $\Hom(x,y) \times_{U(x)} \widetilde{U}(x)$ we have the equality $p_u \circ q(f,u) \circ s = f$ and by Lemma \ref{lemma-pullback-of-diagonal-section} the equality $\widetilde{U}(q(f,u) \circ s) \delta(u) = s$, we conclude $\psi \circ \phi = \id$. Thus, $\phi$ is a bijection satisfying that $(p_u \circ -) \circ \phi$ is the first projection and ${\delta(u)}_x \circ \phi$ is the second projection, showing that our square is a pullback square.   
\endproof

\begin{lemma}\label{lemma-preservation-of-pullback-squares}
The functor $M_!$ maps the distinguished pullback squares based on $\partial$ to pullback squares in $\widehat{\cal C}$.
\end{lemma}
\proof
We need to prove that the image under $M_!$ of a pullback square of the form
\[
\begin{tikzcd}
(P;\eta) \arrow[r, "Q(\eta)"] \arrow[d, "p_{P, \eta}"] & 
\widetilde{U} \arrow[d, "\partial"] \\
P \arrow[r, "\eta"] & 
U
\end{tikzcd}
\]
is a pullback square in $\widehat{\cal C}$. We let the presheaves $P_{(y,u)}$'s be given by the following pullback squares.
\[
\begin{tikzcd}
P_{(y,u)} \arrow[r] \arrow[d] &
y \arrow[d, "c_{(y,u)}"] \\
P \arrow[r, "\eta"] &
U
\end{tikzcd}
\]
Next, we know from Lemma \ref{lemma-pullback-of-univ-along-cocone} that the following square is a pullback square
\[
\begin{tikzcd}[column sep=large, row sep=large]
u \arrow[r, "\delta(u)"] \arrow[d, "Y(p_u)"] & 
\widetilde{U} \arrow[d, "\partial"] \\
y \arrow[r, "c_{(y,u)}"] & 
U,
\end{tikzcd}
\] 
hence we have the following diagram composed of two pullback squares. 
\[
\begin{tikzcd}[column sep=large]
P_{(y,u)} \arrow[r, "{c_{(y,u)}}^* \eta"] \arrow[d, "\eta^* c_{(y,u)}"]& 
y \arrow[d, "c_{(y,u)}"] & 
u \arrow[d, "\delta(u)"] \arrow[l, "Y(p_u)"'] \\
P \arrow[r, "\eta"] & U & \widetilde{U} \arrow[l, "\partial"']
\end{tikzcd}
\]
Now, we write each $P_{(y,u)}$ as a colimit of representables
\[
\varinjlim_{z} P_{yz} \coloneqq \varinjlim_{(z,v)\in \el(P_{(y,u)})^{\op}} z
\]
Since in $\widehat{\cc}$ pulling back commutes with colimits, we have
$$\mld \varinjlim_{y} \varinjlim_{z} (P_{yz} \times_{y} u) \cong & \varinjlim_{y} (P_{(y,u)} \times_{y} u)\\ \cong P\times_U \widetilde{U}$$
and by the same argument, since $M_!$ preserves colimits, we have
\[
\varinjlim_{y} \varinjlim_{z} (M_!(P_{yz}) \times_{M_!(y)} M_!(u)) \cong M_!(P) \times_{M_!(U)} M_!(\widetilde{U}).
\]
So, in order to conclude, $M_!$ being colimit-preserving, it suffices to prove that we have
\[
M_!(P_{yz} \times_{y} u) \cong M_!(P_{yz}) \times_{M_!(y)} M_!(u).
\]
But $M_!(P_{yz} \times_{y} u)$ is the image under $M_!$ of the pullback of $Y_{\cc}(p_u)$ along $Y_{\cc}(f)$, where $f\colon z \rightarrow y$ is the unique morphism of $\cc$ such that $Y_{\cc}(f)$ is the composition 
\[
P_{yz} \rightarrow P_{(y,u)} \xrightarrow{{c_{(y,u)}}^*\eta} y,
\] the Yoneda embedding being fully faithful. This last pullback is isomorphic to the image under $Y_{cc}$ of the distinguished square
\[
\begin{tikzcd}[column sep=large, row sep=large]
f^* u \arrow[d, "p_{f^*u}"]   \arrow[r, "{q(f,u)}"]& 
u \arrow[d, "p_u"] \\
z \arrow[r, "f"] &
y
\end{tikzcd}
\]
in $\cc$. Since $M_! \circ Y_{\cc} \cong Y_{\cal C} \circ M$ (\textit{cf.} \cite[Proposition 3.7.3]{handbookcatalg}), we conclude $M_!(P_{yz} \times_{y} u) \cong M_!(P_{yz}) \times_{M_!(y)} M_!(u)$ using the assumption that $M$ maps the distinguished pullback squares of $\cc$ to pullback squares of $\cal C$ and the fact that $Y_{\cal C}$ preserves pullback squares.    
\endproof

\begin{proposition}\label{proposition-morphism-of-univ-cat}
The triple $(M_!, \id, \id)$ is a morphism of universe categories from $(\widehat{\cc}, \partial)$ to $(\widehat{\cal C}, \partial')$.
\end{proposition}
\proof
It follows from Lemma \ref{lemma-preservation-of-final-object} and Lemma \ref{lemma-preservation-of-pullback-squares}.
\endproof

\subsection{Lifting functors to morphisms of C-systems}\label{subsec-lifting}

\begin{theorem}\label{theorem-main-result}
Let $\cal C$ be a category, $\cc$ be a C-system and $M \colon \cc \sr {\cal C}$ a functor such that $M(\pt)$ is a final object of $\cal C$ and $M$ maps the distinguished pullback squares of $\cc$ to pullback squares of $\cal C$. Then there exists a universe category $(\widehat{\cal C}, p_M)$ and a C-system morphism $M'\colon \cc \sr \cc(\widehat{\cal C}, p_M)$ such that the square
\[
\begin{tikzcd}
\cc \arrow[r, "M"] \arrow[d, "M'"]
& {\cal C} \arrow[d, "Y_{\mathcal{C}}"] \\
\cc(\widehat{{\cal C}}, p_M) \arrow[r, "\int"] 
& \widehat{\mathcal{C}}
\end{tikzcd}
\]
commutes up to a functor isomorphism, with $Y_{\cal C}$ denoting the Yoneda embedding.
\end{theorem}
\proof
Constructions \ref{construction-univ-cat} and \ref{construction-of-generated-C-system} provide a C-system $\cc(\widehat{\mathcal{C}}, \partial')$ and Proposition \ref{proposition-morphism-of-univ-cat} and Construction \ref{construction-lift-of-functor-of-univ-cat} provide a morphism of C-systems $H\coloneqq \cc(M_!, \id, \id)$ from $\cc(\widehat{\cc}, \partial)$ to $\cc(\widehat{\cal C}, \partial')$. Define $M'\colon \cc \sr \cc(\widehat{\mathcal{C}}, \partial')$ as $H \circ I_{\cc}$. Lemma  \ref{lemma-natural-iso-with-functor-of-univ-cat} applied to $F\coloneqq M_!$ provides a natural isomorphism $\psi \colon \int \circ H \sr M_! \circ \int$, while Lemma \ref{lemma-natural-iso-with-yoneda} provides a natural isomorphism $\psi'\colon \int \circ I_{\cc} \sr Y_{\cc}$ and thus we define a natural isomorphism 
\[
\psi'' \colon \int \circ H \circ I_{\cc} \sr M_! \circ Y_{\cc}
\]
with component at $x$ in $\ob(\cc)$ given by the following formula.
\[
\psi''_x \coloneqq M_!(\psi'_x) \circ \psi_{I_{cc}(x)}
\]
Since $M_! \circ Y_{\cc}$ is isomorphic to $Y_{\cal C} \circ M$,
we finally obtain a natural isomorphism from $\int \circ H \circ I_{\cc}$ to $Y_{\mathcal{C}} \circ M$, \textit{i.e} a natural isomorphism from $Y_{\mathcal{C}} \circ M$ to $\int \circ M'$ as required.
\endproof

\end{document}